\documentclass[12pt]{article}
\usepackage{amsmath}
\usepackage{amssymb}
\usepackage{amsthm}
\usepackage{mathrsfs}
\usepackage[normalem]{ulem}
\usepackage[all]{xy}
\usepackage{multicol}
\usepackage{color}
\usepackage{lipsum}
\usepackage{enumitem}
\usepackage{wasysym}
\usepackage[utf8]{inputenc}
\usepackage{longtable} 
 \usepackage{sans}

\theoremstyle{theorem}

\theoremstyle{definition}

\newlength{\bibitemsep}
\newlength{\bibparskip}
\let\oldthebibliography\thebibliography
\renewcommand\thebibliography[1]{%
  \oldthebibliography{#1}%
  \setlength{\parskip}{\bibitemsep}%
  \setlength{\itemsep}{\bibparskip}%
}

\title{\fontsize{14}{14}  \vspace*{-40pt}\textsc{Totalities of Infinite Sets  \vspace{-3ex}}}
\date{\vspace{-5ex}}
\author{ \fontsize{12}{12} \textsc{William Johnston}}%

\definecolor{DarkRed}{cmyk}{0, 0.89, 0.85, 0.35}
\begin{document}
 \maketitle  

\frenchspacing

\begin{abstract} \noindent {Four constructions result from a desire to create enhancements to Cantor's infinite real set cardinality. Each continues to keep Cantor's cardinality formulation in place while providing new comparisons of arbitrary infinite sets. To distinguish their features from a single determinant of a bijection between sets, three of the constructions are characterized by the term ``totality'' in place of cardinality.  We use outer measure in one construction to attain a comparison between subsets of an arbitrary metric space $X$ vis-a-vis subsets of $\mathbb{P}(X)$.}
\end{abstract}

\noindent \section{Introduction}

\noindent Georg Cantor's definition of infinite set cardinalities occured in 1874 with  (translated from German) ``On a Property of the Collection of All Real Algebraic Numbers'' \cite{Cantor}, using the paradigm that two sets have the same cardinality when there exists a bijection between them. Assuming the continuum hypothesis, in the Zermelo-Fraenkel axiomatic system infinite real sets have either countable $\aleph_0$ or uncountable $\aleph_1$ cardinality.\ To offer an enhancement that complements rather than replaces Cantor's definition, we start with a new definition for real set cardinality, based on both the cardinality of the set and the cardinality of the complementary set. Its construction produces more than two different infinite cardinality amounts.\\    

\noindent {\bf Definition 1.}\  A nonempty real set $A$ is countable when $A$ can be written as a sequence $A=\{a_1, a_2, \ldots \}$, either finite (of length, say, $n$) or infinite in length.  
\begin{itemize}
\item When the set $A$ is finite with $n$ elements, we write $|A|=n$, including the case of the empty set $|\phi|=0$.\  When the set is countably infinite, we write $|A|=\aleph_0$.  
\item When $A$ is uncountable, we write, using a new notation, $|A|=\aleph_{1}\!\backslash\aleph_1$, except for the cases where its complement $ \mathbb{R} \backslash A$ is either finite or countably infinite.\ 
\item When the complement is countably infinite we write $|A|=\aleph_{1}\!\backslash\aleph_0$, and when the complement is finite of size $n$, we write $|A|=\aleph_1\!\backslash n$, including the case $|\mathbb{R}|=\aleph_1\!\backslash 0$. 
\end{itemize}
For any $n=0,1,2,3,\ldots$, a natural comparison is \[n < \aleph_0 < \aleph_{1}\!\backslash\aleph_1<\aleph_{1}\!\backslash\aleph_0< \aleph_{1}\!\backslash n. \]
  In addition, two sets that have the same cardinality have a bijection between them {\it and} a bijection between their complements.\\

\noindent \underline{Examples}: 1.\ $|\{7,  \sqrt{2}, \pi, 100 \} | = 4 $.  2.\   $|\mathbb{N} | =|\mathbb{Q} | = \aleph_{0}$.\ 3.\ $|(0,1)|=\aleph_{1}\!\backslash\aleph_1$.\ 4.\ 
$|\mathbb{I}|= \aleph_{1}\!\backslash\aleph_0$ for $\mathbb{I}$ the irrationals.\ 5.\ 
$|\mathbb{R}\backslash \{7,  \sqrt{2}, \pi, 100 \}  | = \aleph_{1}\!\backslash 4$.\  \\ 

An example that promotes an interest in this new definition is its application to the set of normal numbers {\bf N} (in base 10), for which no digit, nor no sequence of digits, occurs more frequently than any other. It is not yet known if $\pi$ is normal, nor $\mathrm{e}$. In 1909 Emil Borel proved in \cite{Borel}, in the sense of Lebesgue measure, that almost all real numbers are normal, and so {\bf N} is uncountable. But, as is commonly known, the bijection $f(0.a_1a_2a_3\ldots)=0.a_10a_20a_30\ldots$ maps $[0,1]$ onto a subset $S$ of non-normal numbers, since, for example, the sequence $rr$ for any nonzero digit $r$ never exists for any number in the output set $S$.  Therefore $\mathbb{R}\backslash${\bf N} is uncountable.  In terms of Definition 1, $|${\bf N}$|= \aleph_{1}\!\backslash\aleph_1$ and the cardinality of normal numbers is strictly smaller than that of irrationals, a statement nicely in tune with how few types of normal numbers have been found.\\ 

\noindent {\bf $\S$ 2.}  \textsc{\bf Using A Step-by-Step Paradigm} \vspace{6pt}

\noindent Imagine imposing paradigms to determine cardinality that go beyond or substitute for the demand that two sets with the same cardinality must have a bijection between them.  In such cases, we use the term ``totality'' of the set in place of ``cardinality'' of the set.  Just as for cardinality, a well-defined totality on the subsets of any set $X$ should:
\begin{itemize}  [itemsep=-2pt]  \vspace{-5pt}
\item be uniquely defined for each subset of $X$;
\item equal the number of elements in the set for any finite set;
\item be monotonic; i.e., for all subsets $A, B \subseteq X$, if  $A \subseteq B$, then the totality of $A$ is less than or equal to the totality of $B$. \vspace{-5 pt}
\end{itemize}
In addition, when $X=\mathbb{R}$, we desire:\begin{itemize}[itemsep=-2pt]  \vspace{-5pt}
\item if $A$ and $B$ have the same totality, then they also have the same cardinality as determined by Cantor.  
\end{itemize}
One useful construction is to employ a step-by-step process to formulate the totality of an infinite set.\ Each step could expand the number of infinite set totalities defined at the immediately previous step, partitioning real-valued sets into more and more divisions of totality.\  Any of the steps could give a well-defined definition of infinite set totalities. In addition, when the steps could continue indefinitely, then an infinite real-valued set totality could be defined as a limit of the step totalities.\  Of course, any finite set totality remains unchanged from its cardinality.\ A simple two-step example follows.\ In this construction and beyond, totalities of infinite sets are labeled using an Omega notation with a subscript size $s$ (as $\Omega_{s}$); this notation distinguishes totality labels from cardinality labels $\aleph_n$.\\

\noindent {\bf Definition 2.\ Infinite real set totality formed from a two-step process.}\  
\begin{itemize}  [itemsep=-2pt]  \vspace{-1pt}
\item In Step 1, define three totalities: $\Omega_0$ if countable; $\Omega_{1/2}$ if uncountable with uncountable complement; and $\Omega_1$ if uncountable with countable complement. 
\item In Step 2, define four totalities, using the Step 1 totalities $\Omega_0$ and $\Omega_1$ but partitioning sets $A$ with Step 1 totality $\Omega_{1/2}$ into two totalities: $\Omega_{1/3}$  when every closed interval $B$ that contains $A$ has $B\backslash A$ uncountable; and $\Omega_{2/3}$ when there exists a closed interval $B$ that contains $A$ with $B\backslash A$ countable.   \vspace{-1pt}
\end{itemize} 
After Step 2, for example the Cantor set $C$ has totality  $|C|=\Omega_{1/3}$, but any nonempty bounded interval $I$ has $|I|=\Omega_{2/3}$.  Note that after each step in Definition 2, each set in a distinct totality still has a bijection between them but additional refinements further partition sets with equal cardinalities into extra totalities.\\

\vspace{0pt}
\noindent {\bf $\S$ 3.}  \textsc{\bf Comparing Totalities of Sets of Different Types.} \vspace{6pt}

\noindent We now use an alternative  paradigm to compare sets of different types, such as comparing subsets of $\mathbb{R}$ with subsets of its power set $\mathbb{P}(\mathbb{R})$. Such comparisons may be considered important because, for example, they may reveal new understandings about how subsets within the two types compare -- beyond the ability to equate their cardinalities via a bijection between them. To compare infinite subsets of any space $X$ against infinite subsets of  $\mathbb{P}(X)$, we desire totality constructions that satisfy four strategic goals:
\begin{enumerate}[itemsep=-2pt] \vspace{-3 pt}
\item The totalities should formulaically follow from only topological properties of $X$. 
\item For any totality value $\Omega_t$ of a subset $\mathcal{P} \subseteq \mathbb{P}(X)$, there should exist a subset $A \subseteq X$ that also has totality $\Omega_t$. Similarly, if $|A|=\Omega_t$ for some $A \subseteq{X}$, then there should exist $\mathcal{P} \subseteq \mathbb{P}(X)$ with $|\mathcal{P}|=\Omega_t$.
\item The totalities of each of the full spaces $X$ and $\mathbb{P}(X)$ should equal. 
\item A totality comparison between any two sets should give a useful meaning about the sets' comparative sizes. 
\end{enumerate}
We call the pair of totalities on $X$ and  $\mathbb{P}(X)$  that satisfy each of these goals a {\it  strategic totality pair}.
This section shows how outer measure can be used to construct one. It thus uses outer measure in a new, broad manner: to formulate totalities that allow totality comparisons between sets of different types in a strategically designed, determined manner 

The required properties of any outer measure is standard (cf. \cite[p. 250]{Royden}). An outer measure $m$ must be defined on all subsets of a set $X$ as an extended real-valued set function $m(A), A \subseteq X$ for which:
\begin{enumerate}[itemsep=-2pt] \vspace{-3 pt}
\item $m(\varnothing) = 0$; 
\item if $A \subseteq B$, then $m(A) \le m(B)$; and 
\item  for subsets $A, B_1, B_2, \ldots \in  \mathbb{P}(X)$, if $ A \subseteq \bigcup_{j=1}^\infty B_j$, then $ m(A) \leq \sum_{j=1}^\infty m(B_j)$. \vspace{-2pt}
\end{enumerate}

To develop comparisons that are uniquely formulated from a strategic totality pair, we restrict attention to sets $X$ that form a topological metric space $(X,d)$ with metric $d$. In such a space, \begin{center}$B_r(x)=\{y \in X\colon d(x,y) < r\}$,\end{center} which are the open balls centered at $x$ with radius $r$,  always form a base for a topology on $X$; i.e., the open sets of $X$ are the unions of open balls. A host of different outer measures can be constructed on any infinite metric space. To construct a uniquely determined outer measure $m$ on subsets of $X$, we begin (cf. \cite[p. 54]{Royden}) by defining the measure of any open ball as its ``diameter length'' so that $m(\, B_r(x)  \, ) = 2r$. Then define $m(A)$ for any $A \subseteq X$ by considering the countable collections $\{ B_n \}$ of open balls that cover $A$; i.e., collections for which $A \subseteq \cup B_n$. For each such collection, set \vspace{-6pt}
$$m(A) = \inf\limits_{A \subseteq \cup B_n} \sum_{n=1}^{\infty} m(B_n) $$
and realize $m$ is an outer measure on $X$ (cf. \cite[Lemma 4 and proof, p. 254]{Royden}).  We label this outer measure the ``{\it metric outer measure}'' on the metric space $(X,d)$. It is uniquely determined by the metric.

The metric outer measure for $X$ naturally produces a set totality for any $A \subseteq X$. When $A$ is finite, its totality is its cardinality. The totality of an infinite countable set $A$ is $|A|=\Omega_{-1}$. The totality of any other infinite set $A$ is $|A|=\Omega_{m(A)}$.  Using the Hausdorff metric, in a similar way the metric $d$ on $X$ also formulaically produces a set totality for the power set  $\mathbb{P}(X)$, which we describe in the next theorem's proof. These constructions allow a strategic totality comparison between any subset of $X$ and any subset of $\mathbb{P}(X)$. \\

\noindent {\bf Theorem 1:} {\it Using the method described above, equip a given metric space $(X,d)$ with its metric outer measure $m$ so that $m(B_r(x))=2r$ and $m(A) = \inf\limits_{A \subseteq \cup B_n} \sum_{n=1}^{\infty} m(B_n) $ for any $A \subseteq X$.   
The Hausdorff metric on $\mathbb{P}(X)$, which is formulaically determined by the metric $d$ on $X$, produces a corresponding metric outer measure on $\mathbb{P}(X)$. The two metric outer measures -- one on $X$ and one  on  $\mathbb{P}(X)$ -- formulate a strategic totality pair. }
\proof  The first part of this proof gives details on how the metric $d$ for $X$ produces a  Hausdorff metric that defines open balls centered at any given set $S \in \mathbb{P}(X)$. For each pair of non-empty subsets $S, T \subseteq X$, the Hausdorff distance between $S$ and $T$ is 
$$ d_{\mathrm H}(S,T) \equiv  \max\left\{\,\sup_{s \in S} d_h(s,T),\ \sup_{t \in T} d_h(S,t) \,\right\},$$ 
 where $d_h(s, T) \equiv \inf\limits_{t \in T} d(s,t)$  quantifies the distance from a point $s \in S$ to the subset  $T$. By convention, $d_h(s, \emptyset)=0$ for any point $s \in X$, and 
$ d_{\mathrm H}(\emptyset,T)=0$. Similarly, $d_h(S, t) \equiv \inf\limits_{s \in S} d(s,t)$  quantifies the distance from a point $t \in T$ to the subset  $S$. By convention, $d_h(\emptyset, t)=0$ for any point $t \in X$.

Any nonempty open ball $B_{r}(S)$ on $\mathbb{P}(X)$ centered at a set $S \in \mathbb{P}(X)$ with Hausdorff radial distance $r$ is then defined as the collection of subsets $T \in \mathbb{P}(X)$ within Hausdorff distance less than $r$, so that:
$$B_{r}(S) \equiv \left\{ T \in \mathbb{P}(X)\colon d_{\mathrm H}(S,T)<r \right\}.$$
Define the ``Hausdorff metric outer measure'' $\mu$ on $\mathbb{P}(X)$ beginning with the definition of $\mu $ on these open balls as $\mu( B_{r}(S) ) \equiv 2r$. Then define the metric outer measure for any subset $\mathcal{P} \subseteq \mathbb{P}(X)$ (where $\mathcal{P} $ is a collection of subsets of $X$) in the standard way:
$$\mu (\mathcal{P}) = \inf \big\{\sum_{n=1}^{\infty} \mu(B_n)\colon   \mathcal{P} \subseteq \bigcup_{n=1}^{\infty} B_n \text{ for } B_n \text{ of the form } B_n= B_r(S) \big\}.$$ Note that the metric $d$ on $X$ determines these metric outer measures for  subsets of $\mathbb{P}(X)$.

We now define the totality of any finite subset to be its cardinality, the totality of an infinite countable set to be $\Omega_{-1}$, the totality of any other infinite subset $\mathcal{P}$ of $\mathbb{P}(X)$ to be $|\mathcal{P}|=\Omega_{\mu(\mathcal{P})}$, and the totality of any other infinite subset $A$ of $X$ to be $|A|=\Omega_{m(A)}$. Each set has a uniquely defined totality because it has a uniquely defined metric outer measure. Well order these infinite set totalities as $\Omega_s < \Omega_t$ exactly when $s<t$. The monotonicity of outer measure forces each of the two  constructions of totality (one for subsets of $X$ and one for subsets of $\mathbb{P}(X)$) to be monotonic; i.e., $\Omega_{s} \le \Omega_{t}$ when $B \subseteq C$ with $|B|=\Omega_{s}$ and $|C|=\Omega_{t}$. Therefore, in each case for $X$ and $\mathbb{P}(X)$, the three requirements for a well-defined totality are satisfied. Furthermore, the facts that  \begin{enumerate}[itemsep=-2pt] \vspace{-3 pt}
\item a meaningful outer measure size comparison via set totalities results between sets of the two vastly different types, and 
\item the open balls in each space generate all totality values from $\Omega_0$ to $\Omega_{m(X)}$, 
\end{enumerate} \noindent confirm the totalities follow the four requirements of a strategic totality pair. \hfill $\blacksquare$\\

\noindent {\bf Example:} Employ Theorem 1 with the absolute value metric on $\mathbb{R}$ to compare the totality of the interval $A=(-1,1)$ in $\mathbb{R}$ with the totality of  $\mathbb{P}({A}) \subseteq \mathbb{P}(\mathbb{R})$. Since $A$ has radius 1, $|A|=\Omega_2$. Using the resulting Hausdorff metric, the open ball centered at $\{ 0 \}$ with radius 1 is $\mathbb{P}(A)$, since each subset of $A$ is strictly within 1 Hausdorf metric distance from $\{ 0 \}$ and any set $T$ not contained in $A$ has $\sup_{t \in T} d_h(\{0\},t) \ge 1$. We conclude $\mu(\mathbb{P}(A))=2$. Hence $A$ and $\mathbb{P}(A)$ have equal totalities $\Omega_2$. \\

\noindent {\bf Commentary:} Given distinct metric spaces $(X,d_1)$ and $(Y,d_2)$, the spaces' metric outer measures and their corresponsing set totalities may allow strategic totality comparisons between infinite subsets taken from $X$ or $Y$.  The situation comparing subsets from $X$ vs. from $\mathbb{P}(X)$ is of special note.  When $X=\mathbb{R}$, real sets with the same totality also have the same cardinality.

Cantor's theorem states the cardinality of $\mathbb{P}(X)$ is always strictly larger than the cardinality of $X$. Totality in Theorem 1 alternatively compares the two spaces; for example, an infinite metric space $X$ having infinite metric outer measure $m(X)$ determines the corresponding Hausdorff metric outer measure $\mu(\, \mathbb{P}(X)\, )$ also as infinity. The two sets have equal totalities $\Omega_{\infty}$. Totality comparison of sets is determined only from the metric $d$ on $X$; every step in the computations of $\Omega_{m(A)}$ and $\Omega_{\mu(B)}$ for arbitrary $A \subset X$ and $B \subset \mathbb{P}(X)$ follows from the metric $d$. As illustrated in the example, Theorem 1 allows the totality of any subset of $X$ to be compared to the uniquely defined totality of any subset of $\mathbb{P}(X)$, which can provide a strong mathematical sense of which sets are large and small in their infinite amount. An analysis of totalities for the power set of the power set of $X$ also follows, as does any finite number of additional iterations of power sets of $X$. Each iteration employs a strategic totality pair. Calculating such comparisons for sets within any metric space can result for $X$ a Hilbert or Banach space such as $L^2(\mathbb{R})$ and its power set. This type of comparison is one item that makes enjoyable the application of Cantor's cardinality to different sets.  Any strategic totality pair, such as in Theorem 1, provides that mathematical enjoyment as well. \\

\noindent {\bf $\S$ 4.} \textsc{\bf Real-Valued Set Totality from an Infinite-Step Process} \vspace{8pt}

This section expands on the totality construction in Section 2, creating a totality for real-valued infinite sets that is defined from a step-by-step process, one that now can be taken out to any number of steps $n$. Each step expands the number of infinite set totalities that were defined in the immediate previous step, partitioning infinite real-valued sets into more and more divisions of totality. Any of the steps gives a well-defined definition of infinite set totalities so that two sets with the same totality also have the same cardinality as determined by Cantor.\ As is standard, a finite set's totality equals its cardinality.\ 

In the description that follows, $|A|_n$ stands for the totality of the infinite real-valued set $A$ at the $n$th step for $n = 1,2, 3, \ldots$.\ We again label totalities of infinite sets as $\Omega_{s}$.\ To simplify totality comparisons, the Omega subscript sometimes uses the letter ``$M$.''\ We start by detailing Step 1 totalities.\vspace{2pt}

\noindent {\bf{Step 1:}} For any countably infinite real-valued set $A$, we write $|A|_1=\Omega_{-1}$, and for any other infinite real-valued set $A$:
\begin{itemize} [itemsep=-2pt]  \vspace{-4pt}
  \item $|A|_1=\Omega_{0}$ when $A$ can be contained in a countable union of open intervals whose total sum of lengths is less than $1/2$.
  \item $|A|_1=\Omega_{1/2}$ when $A$ can be contained\footnote{Of course if $A$ can be contained in a countable union of open intervals whose total sum of lengths is less than $1/2$, then it can be contained in a countable union of open intervals whose total sum of lengths is less than 1.\ The language here should be understood to say that the total sum of lengths fits into the described interval but not a smaller one.} in a countable union of open intervals whose total sum of lengths is at least $1/2$ and less than $1$.
  \item $|A|_1=\Omega_{1}$ when $A$ can be contained in a countable union of open intervals whose total sum of lengths is at least $1$ and less than $3/2$.
  \item $|A|_1=\Omega_{3/2}$ when $A$ can be contained in a countable union of open intervals whose total sum of lengths is at least $3/2$ and less than $2$.
  \item $|A|_1=\Omega_{M-3/2}$ when the complement set $A^c$ can be contained in a countable union of open intervals whose total sum of lengths is at least $3/2$ and less than $2$.
  \item $|A|_1=\Omega_{M-1}$ when $A^c$ can be contained in a countable union of open intervals whose total sum of lengths is at least $1$ and less than $3/2$.
  \item $|A|_1=\Omega_{M-1/2}$ when $A^c$ can be contained in a countable union of open intervals whose total sum of lengths is at least $1/2$ and less than $1$.
  \item $|A|_1=\Omega_{M-0}$ when $A^c$ can be contained in a countable union of open intervals whose total sum of lengths is less than  $1/2$.
  \item $|A|_1=\Omega_{M}$ for all other infinite real-valued sets.\ (This is the ``middle'' infinite totality.)
\end{itemize}
\noindent As examples, the Cantor set $C$ has totality  $|C|_1=\Omega_{0}$.\ Also,  $|[0,1]|_1=\Omega_{1}$, $|(0,\infty)|_1=\Omega_{M}$, $|(-\infty,0) \cup (1,\infty)|_1=\Omega_{M-1}$, and $|\mathbb{R}|_1=\Omega_{M-0}$.\vspace*{6pt}

\noindent Infinite set totality at the general $n$th step works similarly.\vspace{6pt}

\noindent\fbox{%
    \parbox{\textwidth}{%
\noindent {\bf{Step \boldmath $n$ \unboldmath for \boldmath $n=1,2,3,\ldots \ $\unboldmath :}} Any countably infinite real-valued set $A$ has $|A|_n=\Omega_{-1}$.\ For any other infinite real-valued set $A$:
\begin{itemize} [itemsep=-2pt]  \vspace{-9pt}
  \item For $k=0, 1,2,\ldots, 2^{2n}-1$,  $|A|_n=\Omega_{k/2^n}$ when $A$ can be contained in a countable union of open intervals whose total sum of lengths is at least $k/2^n$ and less than $(k+1)/2^n$.\ 
  \item For $k=0,1,2, \ldots, 2^{2n}-1$, $|A|_n=\Omega_{M-k/2^n}$ when the complement set $A^c$ can be contained in a countable union of open intervals whose total sum of lengths is  at least $k/2^n$ and less than $(k+1)/2^n$.\ 
  \item $|A|_n=\Omega_{M}$ for all other infinite real-valued sets.\ (This is the middle infinite cardinality.)\vspace{-9pt}
\end{itemize}
   }%
}
 
\vspace*{8pt} A limit process also determines, as a last offering of this definition’s type, an infinite real-valued set’s totality.\vspace{8pt}

\noindent {\bf Totality as a Limit:} For any infinite real-valued set $A$, when $|A|_n = \Omega_{s_n}$ at each step n, we define \begin{center} $\displaystyle |A|_\infty = \lim_{n \to \infty}|A|_n\equiv \Omega_{s_\infty}$, where $\displaystyle s_\infty = \lim_{n \to \infty} s_n$. \end{center}

After any step and including totalities as a limit, sets with the same totality also have the same cardinality; refinements partition equal cardinality sets into additional totalities. Though different, totality defined as a limit is related to Lebesgue outer measure and always exists. At each step and in the limit, these uncountable totalities are of three forms: \vspace{-0pt}
\begin{center} $\Omega_a, \ \   \ \Omega_M, \   \mbox{ and } \ \  \Omega_{M-a}$, where $a$ is nonnegative finite.\vspace{-0pt} \end{center} 
A comparison of these totalities is also defined: For nonnegative reals $a<b$: \begin{center} $\Omega_a < \Omega_b$, \ \ $\Omega_{M-b} < \Omega_{M-a}$, \ \  $\Omega_a < \Omega_M$, \  and \ $\Omega_M < \Omega_{M-a}$.\end{center}

\noindent The following theorem results.\\

\noindent{\bf  Theorem 2:} {\it For the Totality as a Limit, no set has a totality between $\Omega_{-1}$ and $\Omega_0$. Between any two uncountable infinite set totalities always lies another.}\\

\noindent {\bf $\S$ 5.} \textsc{\bf Conclusion} \vspace{8pt}

\noindent A thought-provoking observation:\ Every consistent axiomatic system has undecidables.\ Using the Totality as a Limit definition in Section 4 may imply the First Incompleteness Theorem has flexibility, as that totality enhances infinite real set cardinality -- an example of a definition refinement  from within an axiomatic system where at the refinement level the continuum hypothesis undecidable is not a concern.\ The continuum hypothesis undecidable is still there, as we continue to use and celebrate Cantor's real-valued infinite set cardinality definition, but the undecidable resolves for the refinement. Though undecidables are always present, might you always be able to refine a definition (not the system) one step at a time to form a resolution at the refinement level for each undecidable the system produces?  Whenever a refinement works in this way, then deductive mathematics becomes more muscular, creating subsidence for such undecidables each time.

\end{document}